\documentclass[12pt]{amsart}
 \usepackage{graphicx}
\usepackage[all]{xypic}
\usepackage{amsmath,graphics}
\usepackage{amsfonts,amssymb}

\theoremstyle{plain}
\newtheorem*{theorem*}{Theorem}
\newtheorem*{lemma*} {Lemma}
\newtheorem*{corollary*} {Corollary}
\newtheorem*{proposition*} {Proposition}
\newtheorem{theorem}{Theorem}[section]
\newtheorem{lemma}[theorem]{Lemma}
\newtheorem{corollary}[theorem]{Corollary}
\newtheorem{proposition}[theorem]{Proposition}
\newtheorem{conjecture}[theorem]{Conjecture}

\theoremstyle{remark}
\newtheorem*{remark}{Remark}
\newtheorem*{definition}{Definition}

\newtheorem*{claim}{Claim}

\theoremstyle{definition}

\def\norm{|| _- ||}

\def\gl{\mbox{GL}}
\def\Q{\Bbb{Q}}
\def\R{\Bbb{R}}
\def\k{\Bbb{K}}
\def\K{\Bbb{K}}

\def\id{\mbox{id}}

\def\Z{\Bbb{Z}}
\def\C{\Bbb{C}}

\def\N{\Bbb{N}}

\def\part{\partial}

\def\a{\alpha}

\def\g{\gamma}

\def\bp{\begin{pmatrix}}

\def\sm{\setminus}
\def\ep{\end{pmatrix}}
\def\bn{\begin{enumerate}}

\def\en{\end{enumerate}}
\def\ba{\begin{array}}
\def\ea{\end{array}}

\def\v{\alpha}
\def\b{\beta}
\def\ti{\tilde}

\def\fr12{\frac{1}{2}}

\def\im{\mbox{Im}}

\def\ker{\mbox{Ker}}

\def\hom{\mbox{Hom}}

\def\deg{\mbox{deg}}

\def\v{\varphi}

\def\upm{[u^{\pm 1}]}

\def\be{\begin{equation}}
\def\ee{\end{equation}}
\def\spm{[s^{\pm 1}]}

\def\gen{\mbox{genus}}

\def\ktbig{\K[t_1^{\pm 1},\dots,t_m^{\pm 1}]}
\def\ftbig{\F[t_1^{\pm 1},\dots,t_m^{\pm 1}]}

\def\ktbigfield{\K(t_1,\dots,t_m)}
\def\kdet{\ktbigfield^\times_{ab}}
\def\tbigpm{[t_1^{\pm 1},\dots,t_m^{\pm 1}]}
\def\tbigfield{(t_1,\dots,t_m)}
\def\F{\Bbb{F}}
\def\k{\Bbb{K}}
\def\K{\Bbb{K}}

\def\kas{\k[s^{\pm 1}]}
\def\ras{\mathcal{R}[s^{\pm 1}]}
\def\ksfield{\k(t)}

\def\ksfield{\k(s)}

\def\ol{}    

\def\K{\Bbb{K}}

\def\gen{\mbox{genus}}

\def\cmtbf#1{} \def\cmt#1{}

\begin{document}

\title[Non--commutative multivariable Reidemeister torsion]
{Non--commutative multivariable Reidemeister torsion and the Thurston norm}
\author{Stefan Friedl and Shelly Harvey}
\date{\today} \address{Rice University, Houston, Texas, 77005-1892}
\email{friedl@alumni.brandeis.edu, shelly@math.rice.edu}
\def\subjclassname{\textup{2000} Mathematics Subject Classification}
\expandafter\let\csname subjclassname@1991\endcsname=\subjclassname \expandafter\let\csname
subjclassname@2000\endcsname=\subjclassname 
\date{\today}
\begin{abstract}
Given a $3$--manifold the second author defined functions
$\delta_n:H^1(M;\Z)\to \N$, generalizing McMullen's Alexander norm,
which give lower bounds on the Thurston
norm. We reformulate these invariants in terms of Reidemeister
torsion over a non--commutative multivariable Laurent
polynomial ring. This allows us to show that these functions are
semi-norms.
 \end{abstract}
\maketitle




\section{Introduction}
Let $M$ be a 3-manifold. Throughout the paper we will assume that all 3-manifolds are compact,
connected and orientable. Let $\phi\in H^1(M;\Z)$.  The \emph{Thurston norm} of $\phi$ is
defined as
 \[
||\phi||_{T}=\min \{ \chi_-(S)\, | \, S \subset M \mbox{ properly embedded surface dual to }\phi\}
\] where given a surface $S$ with connected components $S_1,\dots,S_k$ we
write $\chi_-(S)=\sum_{i=1}^k\max\{0,-\chi(S_i)\}$. We refer to
 \cite{Th86} for details.


Generalizing work of Cochran \cite{Co04} the second author introduced in \cite{Ha05} a function
\[ \delta_n:H^1(M;\Z)\to \N_0\cup \{-\infty\}\]
for every $n\in \N$ and showed that $\delta_{n}$ gives a lower bound on the
Thurston norm for every $n$.  These functions are invariants of the 3-manifold and generalize the Alexander norm defined by C. McMullen in \cite{Mc02}.  We point out that the definition we
use here differs slightly from the original definition when $n=0$
and a few other special cases. We refer to Section
\ref{section:rhon} for details.

The relationship between the functions
$\delta_n$
and the Thurston norm was further strengthened in
\cite{Ha06} (cf. also \cite{Co04} and \cite{Fr05}) where it was shown that the
 $\delta_n$
 give a never decreasing series of lower bounds on the Thurston norm, i.e. for any $\phi\in H^1(M;\Z)$ we
have
\[  \delta_0(\phi) \leq  \delta_1(\phi) \leq \delta_2(\phi)\leq \dots \leq ||\phi||_T.\]
Furthermore it was shown in \cite{FK05c} that under a mild assumption these inequalities are an
equality modulo 2.

Thurston \cite{Th86} showed in particular that $||-||_T$  is a
seminorm. It is therefore a natural question to ask whether the
invariants $\delta_n$ are seminorms as well. In \cite{Ha05} this
was shown to be the case for $n=0$. The following theorem, which is
a special case of the main theorem of this paper (cf. Theorem
\ref{thm:harveynorm}), gives an affirmative answer for all $n$.

\begin{theorem} \label{thm:mainintro}
Let $M$ be a 3--manifold with empty or toroidal boundary. Assume
that
$\delta_n(\phi)\ne -\infty$ for some $\phi\in H^1(M;\Z)$, then
\[ \delta_n:H^1(M;\Z)\to \N_0\]
is a seminorm.
\end{theorem}

This in particular allows us to show that the sequence $\{\delta_{n}\}$ is eventually constant.  That is, there exists an $N\in \N$ such that $\delta_{n}=\delta_{N}$ for all $n\geq N$
(cf. Proposition \ref{prop:eventconstant}).

Initially we discuss a more algebraic problem. Recall that given a
multivariable  Laurent polynomial ring $\ftbig$ over a commutative
field $\F$ we can associate to any non--zero $f=\sum_{\a\in
\Z^m}a_\a t^{\a} \in \ftbig$ a seminorm on $\hom(\Z^m,\R)$ by
\[ ||\phi||_f:=\sup\{ \phi(\a)-\phi(\b) | a_\a\ne 0, a_\b\ne 0\}.\]
Furthermore, to any square matrix $B$ over $\ftbig$ with $\det(B)\neq 0$ we can associate a norm using $\det(B)\in
\ftbig$.

Generalizing this idea to the non-commutative case, in Section \ref{section:multivar} we introduce the notion of a
\emph{multivariable skew Laurent polynomial ring} $\ktbig$ of rank
$m$ over a skew field $\K$. Given a square matrix $B$ over $\ktbig$
we can study its Dieudonn\'e determinant $\det(B)$ which is an
element in the abelianization of the multiplicative group
$\ktbigfield\sm \{0\}$ where $\ktbigfield$ denotes the quotient
field of $\ktbig$. This determinant will in general not be
represented by an element in $\ktbig$. Our main technical result
(Theorem \ref{thm:kstdet}) is that  nonetheless  there is a natural
way to associate a norm to $B$ which generalizes the commutative
case.

Given a $3$--manifold $M$ and  a `compatible'--representation 
$$\pi_1(M)\to \gl(\ktbig,d)$$ we will show in Section~\ref{section:rt}
that the corresponding  Reidemeister torsion can be viewed as a
matrix over $\ktbig$. We will show in Section~\ref{section:rhon}
that for appropriate representations the norm which we can associate
to the matrix over $\ktbig$ agrees with
$\delta_n$. In particular, this implies Theorem
\ref{thm:mainintro}.
We conclude this paper with examples of links for which we compute
the Thurston norm using these invariants.

As a final remark we point out that the results in this paper
completely generalize the results in \cite{FK05b}. Furthermore the
results can easily be extended to studying  2--complexes together
with the Turaev norm which is  modeled on the definition of the
Thurston norm of a 3--manifold. We
refer to \cite{Tu02a} for details.\\

{\bf Acknowledgments:} The  authors would like to thank Tim Cochran,
John Hempel, Taehee Kim and  Chris Rasmussen  for helpful
conversations.

\section{The non--commutative Alexander norm}
%

\subsection{Multivariable Laurent polynomials}\label{section:multivar}
Let $\mathcal{R}$ be a (non--commutative) domain and
$\g:\mathcal{R}\to \mathcal{R}$ a ring homomorphism. Then we denote
by $\ras$ the \emph{one--variable skew Laurent polynomial ring over
$\mathcal{R}$}. Specifically  the elements in $\ras$ are formal sums
$\sum_{i=m}^n a_is^i$ ($m\leq n\in \Z$) with $a_i\in \mathcal{R}$.
Addition is given by addition of the coefficients, and
multiplication is defined using the rule $s^ia=\g^i(a)s^i$ for any
$a\in \mathcal{R}$ (where $\g^i(a)$ stands for $(\g\circ \dots \circ
\g)(a)$). We point out that any element $\sum_{i=m}^n
a_is^i\in \ras$ can also be written uniquely in the form
$\sum_{i=m}^n s^i\ti{a}_i$, indeed, $\ti{a}_i=s^{-i}a_is^i\in
\mathcal{R}$.

In the following let $\K$ be a skew field. We then define \emph{multivariable skew Laurent
polynomial ring of rank $m$ over $\K$} (in non--commuting variables)
to be a ring $R$ which is an algebra over $\K$ with
unit (i.e. we can view $\K$ as a subring of $R$) together with a
decomposition $R=\oplus_{\a\in \Z^m} V_\a$ such that the following
hold: \bn
\item $V_\a$ is a one--dimensional $\K$--vector space,
\item $V_\a \cdot V_\b=V_{\a+\b}$,
\item $V_{(0,\dots,0)}=\K$.
\en In particular $R$ is $\Z^m$--graded. Note that these properties
imply that any $V_\a$ is invariant under left and right
multiplication by $\K$, that any element in $V_\a\sm \{0\}$ is a
unit, and that $R$ is a (non--commutative) domain.

The example to keep in mind is a commutative Laurent polynomial ring $\F[t_1^{\pm
1},\dots,t_m^{\pm 1}]$. Let $t^{\a}:=t_1^{\a_1}\cdot \dots \cdot t_m^{\a_m}$ for
$\a=(\a_1,\dots,\a_m)$, then $V_\a=\F t^\a, \a\in\Z^m$ has the required properties.

Let $R$ be a multivariable skew Laurent polynomial ring of rank $m$
over $\K$. To make our subsequent definitions and arguments easier
to digest we will always pick $t^{\a}\in V_\a\sm \{0\}$ for $\a\in
\Z^m$. It is easy to see that we can in fact pick $t^{\a}, \a\in
\Z^m$ such that $t^{n\a}=(t^{\a})^{n}$ for all $\a \in \Z^m$ and
$n\in \Z$. Note that this choice in particular implies that
$t^{(0,\dots,0)}=1$. We get the following properties: \bn
\item $t^{\a}t^{\ti{\a}}t^{-(\a+\ti{\a})}\in \K^\times$ for all $\a,\ti{\a}\in \Z^m$, and
\item $t^\a \K =\K t^{\a}$ for all $\a$.
\en This shows that the notion of multivariable skew Laurent
polynomial ring of rank $m$ is a generalization of the notion of
twisted group ring of $\Z^m$ as defined in \cite[p.~13]{Pa85}. If
$m=1$ then we have $t^{(n)}\in V_{(n)}$ such that
$t^{(n)}=(t^{(1)})^n$ for any $n\in \Z$. We write $t^n=t^{(n)}$. In
particular we have a one--variable skew Laurent polynomial ring as
above.

The argument of \cite[Corollary~6.3]{DLMSY03} can be used to show that any such Laurent polynomial
ring
 is a (left and right) Ore domain and in particular has a (skew) quotient field.
We normally denote a multivariable skew Laurent polynomial ring of rank $m$ over $\K$ suggestively
by $\K[t_1^{\pm 1},\dots,t_m^{\pm 1}]$ and we denote the quotient field of  $\K[t_1^{\pm
1},\dots,t_m^{\pm 1}]$ by $\K(t_1,\dots,t_m)$.

\subsection{The Dieudonn\'e determinant}

In this section we recall several well--known definitions and facts.
Let $\mathcal{K}$ be a skew field. In our applications $\mathcal{K}$
will be the quotient field of a multivariable skew Laurent
polynomial ring. First define
$\gl(\mathcal{K}):=\underset{\rightarrow}{\lim} \,
\gl(\mathcal{K},n)$, where we have the following maps in the direct
system: $\gl(\mathcal{K},n)\to \gl(\mathcal{K},n+1)$ given by $A
\mapsto \bp A &0\\0&1\ep $,
then define
$K_1(\mathcal{K})=\gl(\mathcal{K})/[\gl(\mathcal{K}),\gl(\mathcal{K})]$.
For details we refer to \cite{Mi66} or \cite{Tu01}.

Let  $A$ a square matrix over $\mathcal{K}$. After elementary row
operations and destabilization we can arrange that in
$K_1(\mathcal{K})$ the matrix $A$ is represented by a $1\times
1$--matrix $(d)$. Then the Dieudonn\'e determinant $\det(A)\in
\mathcal{K}_{ab}^\times:=\mathcal{K}^\times/[\mathcal{K}^\times,\mathcal{K}^\times]$
(where $\mathcal{K}^\times:=\mathcal{K}\sm \{0\}$) is defined to be
$d$.
It is well--known that the Dieudonn\'e determinant induces an
isomorphism $\det:K_1(\mathcal{K})\to \mathcal{K}_{ab}^\times$. We
refer to \cite[Theorem~2.2.5 and Corollary~2.2.6]{Ro94} for more
details.

\subsection{Multivariable skew Laurent polynomial rings and
seminorms}\label{section:polynorm} Let $\kas$ be a
one--variable skew Laurent polynomial ring and let $f\in \kas$. If
$f=0$ then we write $\deg(f)=-\infty$, otherwise, for
$f=\sum_{i=m}^n a_is^i\in \kas$ with $a_m\ne 0, a_n \ne 0$ we define
$\deg(f):=n-m$. This extends to a homomorphism $\deg:\ksfield\sm
\{0\}\to \Z$ via $\deg(fg^{-1})=\deg(f)-\deg(g)$. Since $\deg$ is a
homomorphism to an abelian group this induces a homomorphism
$\deg:\ksfield_{ab}^\times \to \Z$. Note that throughout this paper
we will apply the convention that $-\infty < a$ for any $a\in \Z$.


For the remainder of this section let $\ktbig$ be a multivariable skew Laurent polynomial ring
of rank $m$ together with a choice of $t^{\a}, \a\in \Z^m$ as above.
 Let $f \in \ktbig$. We can
write $f = \sum_{\a\in \Z^m} a_{\a}t^\a$ for some $a_{\a}\in \k$. We associate a seminorm $\norm_f$
on $\hom(\R^m,\R)$ to $f$ as follows. If $f=0$, then we set $\norm_f:=0$. Otherwise we set
\[  ||\phi||_{f} :=\sup \{ \phi(\a)-\phi(\b) | a_\a\ne 0, a_\b\ne 0 \}. \]
Clearly $\norm_f$ is a seminorm and does not depend on the choice of $t^\a$. This seminorm
should be viewed as a generalization of the degree function.

Now let $\tau \in K_1(\ktbigfield)$ and let $f_n,f_d \in \ktbig\sm \{0\}$ such that
$\det(\tau)=f_nf_d^{-1} \in \ktbigfield^\times_{ab}$. Then define
\[ ||\phi||_{\tau}:=\max\{0,||\phi||_{f_n}-||\phi||_{f_d}\} \]
for any $\phi \in \hom(\R^m,\R)$. By the following proposition  this
function is well--defined.

\begin{proposition}\label{prop:welldef}
Let $\tau \in K_1(\ktbigfield)$. Let $f_n,f_d, g_n,g_d \in \ktbig\sm \{0\}$ such that
$\det(\tau)=f_nf_d^{-1}=g_ng_d^{-1} \in \ktbigfield^\times_{ab}$. Then
\[ \norm_{f_n}-\norm_{f_d}=\norm_{g_n}-\norm_{g_d}.  \]
\end{proposition}

We postpone the proof  to Section \ref{section:proofwelldef}.

Let $B$ be a matrix defined over $\ktbig$. Then it is in general not
the case that $\det(\ktbig)$ can be represented by an element in
$\ktbig$. But we still have the following result which is the main
technical result of this paper.

\begin{theorem}\label{thm:kstdet}
If $\tau \in K_1(\ktbigfield)$  can be represented by a matrix defined over $\ktbig$, then
$\norm_\tau$  defines a seminorm on $\hom(\R^m,\R)$.
\end{theorem}

We postpone the proof  to Section \ref{section:proofkstdet}.\\

Now let $\phi:\Z^m\to \Z$ be a non--trivial homomorphism. We will
show that  $||\phi||_B$ can also be viewed as the degree of a
polynomial associated to $B$ and $\phi$. We begin with some
definitions.
 Consider
\[ \K[\ker(\phi)]:=\bigoplus_{\a\in Ker(\phi)} \K t^\a \subset \ktbig. \]
This clearly defines a subring of $\ktbig$ and the argument of \cite[Corollary~6.3]{DLMSY03}
shows that $\K[\ker(\phi)]$ is an Ore domain with skew field which we denote by
$\K(\ker(\phi))$.

Let $d\in \Z$ such that $\im(\phi)=d\Z$ and pick $\b=(\b_1,\dots,\b_m)\in \Z^m$ such that
$\phi(\b)=d$.
Let $\mu:=t^\b$. Then we can form one--variable Laurent polynomial rings
$(\K[\ker(\phi)])[s^{\pm 1}]$ and $\K(\ker(\phi))[s^{\pm 1}]$ where $sk:=\mu k \mu^{-1}\, s$ for
all $k\in \K[\ker(\phi)]$ respectively for all $k\in \K(\ker(\phi))$. We get a map
\[ \ba{rcl} \g_\phi: \ktbig &\xrightarrow{\cong}& (\K[\ker(\phi)])[s^{\pm 1}]\\
  \sum_{\a\in \Z^m} k_\a t^{\a}&\mapsto&\sum_{\a\in \Z^m}k_\a t^{\a}\mu^{-\phi(\a)/d}\, s^{\phi(\a)/d},\ea
\]
 where $k_\a\in \K$ for all $\a\in \Z^m$. Note that
$k_\a t^{\a}\mu^{-\phi(\a)/d}\in \k[\ker(\phi)]$. An easy
computation shows that $\g_\phi$ is an isomorphism of rings. Clearly
we also get an induced isomorphism $\ktbigfield \xrightarrow{\cong}
(\K(\ker(\phi)))(s)$.

Let $B$ a matrix over $\ktbigfield$. Define
$\deg_{\phi}(B):=\deg(\det(\g_\phi(B)))$ where we view $\g(B)$ as a
matrix over $\K(\ker(\phi))(s)$.

\begin{theorem}\label{thm:kstdet2}
Let $B$ a matrix over $\ktbigfield$. Let $\phi \in \hom(\Z^m,\Z)$ non--trivial and  let $d\in
\N$ such that $\im(\phi)=d\Z$. Then
\[ ||\phi||_{B}=d\, \max\{0,\deg_{\phi}(B)\}.  \]
\end{theorem}

Note that this shows in particular that $\deg_{\phi}(B)$ is independent of the choice of $\b$. This
theorem is a generalization of \cite[Proposition~5.12]{Ha05} to the non--commutative case.

\begin{proof}
Since $\g$ and $\deg$ are homomorphisms it is clearly enough to show
that for any $g\in \k\tbigpm\sm \{0\}$ we have
\[ ||\phi||_{g}=d\, \deg(\g_\phi(g)).\]
Write $g=\sum_{\a\in \Z^m}a_{\a} t^{\a}$ with $a_\a\in \K$. Let $d,\b,\mu$ and  $\g: \ktbig
\xrightarrow{\cong} (\K[\ker(\phi)])[s^{\pm 1}]$ as above. Note that $\ker(\phi)\oplus \Z \b
=\Z^m$, hence
\[ \ba{rcl}
g &= &\sum_{i\in \Z}  \sum_{\a\in Ker(\phi)} a_{\a+i\b}t^{\a+i\b},\\
\g_\phi(g) &= &\sum_{i\in \Z} \Big( \sum_{\a\in Ker(\phi)}
a_{\a+i\b}t^{\a+i\b}\mu^{-i}\Big) s^{i}.\ea
\]
 Note that $a_{\a+i\b}t^{\a+i\b}\mu^{-i}\subset
\K t^{\a}$. Since $\k[\ker(\phi)]=\oplus_{\a\in Ker(\phi)}\K t^{\a}$ we get the following
equivalences:
\[ \ba{rrcll} &\sum_{\a\in Ker(\phi)} a_{\a+i\b}t^{\a+i\b}\mu^{-i}&=& 0\\
\Leftrightarrow& a_{\a+i\b}t^{\a+i\b}\mu^{-i}&=&0& \mbox{for all $\a\in \ker(\phi)$}\\
\Leftrightarrow&  a_{\a+i\b}&=&0& \mbox{for all $\a\in \ker(\phi)$}.\ea \]
 Therefore
\[ \ba{rcll} ||\phi||_g&=&& d\, \max_{i\in \Z} \{ \mbox{there exists }\a \in \ker(\phi)\mbox{ such that
}a_{\a+i\b}\ne 0 \} \\
&&-& d\,\min_{i\in \Z} \{ \mbox{there exists }\a \in \ker(\phi) \mbox{ such that
}a_{\a+i\b}\ne 0 \} \\
&=&&  d\,\max_{i\in \Z} \{\sum_{\a\in ker(\phi)}  a_{\a+i\b}t^{\a+i\b}\mu^{-i} \ne 0 \} \\
&&-& d\,\min_{i\in \Z} \{\sum_{\a\in ker(\phi)} a_{\a+i\b}t^{\a+i\b}\mu^{-i} \ne 0 \} \\
&=&&d\,\deg(\g_\phi(g)). \ea \]
\end{proof}

\subsection{Proof of Proposition \ref{prop:welldef}} \label{section:proofwelldef}
We start out with the following basic lemma.
\begin{lemma}\label{lem:newton}
Let $f,g\in \ktbig\sm \{0\}$, then $\norm_{fg}=\norm_{f}+\norm_g$.
\end{lemma}

This lemma is well--known. It follows from the fact that the Newton
polytope of non--commutative multivariable polnyomials $fg$ is the
Minkowski sum of the Newton polytopes of $f$ and $g$.

%

\begin{lemma}
Let $d\in \ktbigfield$ and let  $f_n,f_d,g_n,g_d \in \ktbig$ such
that $d=f_nf_d^{-1}=g_ng_d^{-1} \in \ktbigfield$. Then \[
\norm_{f_n}-\norm_{f_d} =\norm_{g_n}-\norm_{g_d}.\] In particular
\[ \norm_d:=\norm_{f_n}-\norm_{f_d}\]
is well--defined.
\end{lemma}

\begin{proof}
 Recall that by the definition of the Ore localization $f_nf_d^{-1}=g_ng_d^{-1} \in
\ktbigfield$ is equivalent to the existence of $u,v\in \ktbig\sm
\{0\}$ such that $f_nu=g_nv$ and $f_du=g_dv$. The lemma  now follows
immediately from Lemma \ref{lem:newton}.
\end{proof}

\begin{lemma}\label{lem:normmult}
Let $d,e\in \ktbigfield$, then
\[ \norm_{de}=\norm_{d}+\norm_{e}.\]
\end{lemma}

\begin{proof}
Pick  $f_n,f_d,g_n,g_d \in \ktbig$ such that $f_nf_d^{-1}=d$ and
$g_ng_d^{-1}=e$. By the Ore property there exist $u,v \in \ktbig\sm
\{0\}$ such that $g_nu=f_dv$. It follows that
\[ f_nf_d^{-1}g_ng_d^{-1}=f_nvu^{-1}g_d^{-1}=(f_nv)(g_du)^{-1}.\]
The lemma now follows immediately from Lemma \ref{lem:newton}.
\end{proof}

We can now give the proof of Proposition \ref{prop:welldef}.

\begin{proof}[Proof of Proposition \ref{prop:welldef}]
Let $B$ be a matrix defining an element $K_1(\ktbigfield)$. Assume
that we have $f_n,f_d,g_n,g_d \in \ktbig$ such that
$\det(B)=f_nf_d^{-1}=g_ng_d^{-1} \in \kdet$. We can lift the
equality $f_nf_d^{-1}=g_ng_d^{-1} \in \kdet$ to an equality
 \be \label{kstequality} f_n f_d^{-1}=\prod_{i=1}^r [a_i,b_i]\,
 g_n g_d^{-1}  \in \ktbigfield^\times \ee
 for some $a_i,b_i\in \ktbigfield$.
It follows from Lemma \ref{lem:normmult} that $\norm_{[a_i,b_i]}=0$.
It then follows from Lemma \ref{lem:normmult} that
$\norm_{f_nf_d^{-1}}=\norm_{g_ng_d^{-1}}$.
\end{proof}

\subsection{Proof of Theorem \ref{thm:kstdet}} \label{section:proofkstdet}
Now let $\tau \in K_1(\ktbigfield)$ which can be represented by a matrix $B$ defined over
$\ktbig$. We will show that $\norm_\tau=\norm_B$  defines a seminorm on $\hom(\R^m,\R)$.

Because of the continuity and the $\N$--linearity of $\norm_B$ it is enough to show that for any
two non--trivial homomorphisms $\phi,\ti{\phi}:\Z^m\to \Z$ we have
\[ ||\phi+\ti{\phi}||_B \leq ||\phi||_B+ ||\ti{\phi}||_B.\]

Let $\phi,\ti{\phi}:\Z^m\to \Z$ be non--trivial homomorphisms.
Let $d\in \Z$ such that $\im(\phi)=d\Z$ and pick $\b$ with
$\phi(\b)=d$. We write $\mu=t^\b$. As in Section
\ref{section:polynorm} we can form $\K[\ker(\phi)]$ and we also have
an isomorphism $ \g_{\phi}: \ktbig \xrightarrow{\cong}
(\K[\ker(\phi)])[s^{\pm 1}]$.

 Consider $\g_{\phi}(B)$, it is defined over
the PID $\K(\ker(\phi))[s^{\pm 1}]$. Therefore we can use elementary
row operations to turn $\g_{\phi}(B)$ into a diagonal matrix with
entries in $\K(\ker(\phi))[s^{\pm 1}]$. In particular we can find
$a_i,b_i\in \K[\ker(\phi)]$ such that
\[ \det(\g_{\phi}(B))=\sum_{i=r_1}^{r_2} s^i a_ib_i^{-1} \]
 Since $\K[\ker(\phi)]$ is an Ore domain we can in fact find a
 common denominator for $a_ib_i^{-1}, i=r_1,\dots,r_2$. More
 precisely, we can find
$c_{r_1},\dots,c_{r_2}\in \K[\ker(\phi)]$ and $d\in\K[\ker(\phi)]$
such that $a_ib_i^{-1}=c_id^{-1}$ for $i=r_1,\dots,r_2$. Now let
$c=\sum_{i=r_1}^{r_2}s^ic_i$. Then
\[ \det(\g_{\phi}(B))=cd^{-1} \in\K(\ker(\phi))(s)^\times_{ab}\]
where
 $c\in \K[\ker(\phi)][s^{\pm 1}]$ and $d\in
\K[\ker(\phi)]$.  Now let $f=\g_{\phi}^{-1}(c)\in \ktbig,
g=\g_{\phi}^{-1}(d)\in \K[\ker(\phi)]$. Then $\det(B)=fg^{-1}$ and
by Proposition \ref{prop:welldef} we have
\[ \norm_B=\norm_{f}-\norm_{g}. \]
The crucial observation  is that $||\phi||_{g}=0$ and $||\phi+\ti{\phi}||_g=||\ti{\phi}||_g$
since $g\in \K[\ker(\phi)]$. It therefore now follows that
\[ \ba{rcl} ||\phi+\ti{\phi}||_B&=&||\phi+\ti{\phi}||_{f}-||\phi+\ti{\phi}||_{g}\\
&=&||\phi+\ti{\phi}||_{f}-||\ti{\phi}||_{g}\\
&\leq&||\phi||_f+||\ti{\phi}||_{f}-||\ti{\phi}||_{g}\\
&=&(||\phi||_f-||\phi||_{g})+(||\ti{\phi}||_{f}-||\ti{\phi}||_{g})\\
&=&||\phi||_B+||\ti{\phi}||_{B}.\ea \] This concludes the proof of
Theorem \ref{thm:kstdet}.

\section{Applications to the Thurston norm} \label{section:rt}

\subsection{Reidemeister torsion} \label{subsection:rt}

Let $X$ be a finite connected  CW--complex. Denote the universal
cover of $X$ by $\ti{X}$. We view $C_*(\ti{X})$ as a right
$\Z[\pi_1(X)]$--module via deck transformations. Let $R$ be a ring.
Let $\v:\pi_1(X)\to \gl(R,d)$ be a representation, this equips $R^d$
with a left $\Z[\pi_1(X)]$--module structure. We can therefore
consider the right $R$--module chain complex
$C_*^\v(X;R^d):=C_*(\ti{X})\otimes_{\Z[\pi_1(X)]} R^d$. We denote
its homology by $H_i^\v(X;R^d)$. If $H_*^\v(X;R^d)\ne 0$, then we
write $\tau(X,\v):=0$. Otherwise we can define the Reidemeister
torsion $\tau(X,\v)\in K_1(R)/\pm \v(\pi_1(X))$. If the homomorphism
$\v$ is clear we also write $\tau(X,R^d)$.

Let $M$ be a manifold. Since Reidemeister torsion only depends on the homeomorphism type of the
space we can
 define $\tau(M,\v)$  by picking any
CW--structure for $M$.  We refer to the excellent book of Turaev \cite{Tu01} for filling in the
details.

\subsection{Compatible homomorphisms and the higher order Alexander norm} In the following let $M$ be a
3--manifold with empty or toroidal boundary,
let $\psi:H_1(M)\to \Z^m$ be an epimorphism, and let $\K[t_1^{\pm
1},\dots,t_m^{\pm 1}]$ be a multivariable skew Laurent polynomial
ring of rank $m$ as in Section \ref{section:multivar}.

A representation $\v:\pi_1(M) \to \gl(\ktbig,d)$ is called \emph{$\psi$--compatible} if for any
$g\in \pi_1(X)$ we have $\v(g)=At^{\psi(g)}$ for some $A\in \gl(\k,d)$. This generalizes
definitions in \cite{Tu02b} and \cite{Fr05}. We denote the induced representation $\pi_1(M)\to
\gl(\ktbigfield,d)$ by $\v$ as well and we consider the corresponding Reidemeister torsion
$\tau(M,\v)\in K_1(\ktbigfield)/\pm \v(\pi_1(M))\cup \{0\}$.

 We say $\v$ is a \emph{commutative representation}
 if there exists a commutative subfield $\F$ of $\K$ such that for all $g$ we have
$\v(g)=At^{\psi(g)}$ with $A$ defined over $\F$ and if $t^\a, t^{\ti{\a}}$ commute for any
$\a,\ti{\a}\in \Z^m$.

\begin{theorem}\label{mainthmnorm}
Let $M$ be a 3--manifold with empty or toroidal boundary.
Let $\psi:H_1(M)\to \Z^m$ be an epimorphism.  Let $\v:\pi_1(M) \to \gl(\ktbig,d)$ be a
$\psi$--compatible representation such that $\tau(M,\v)\ne 0$. If one of the following holds:
\bn
\item $\v$ is commutative,
\item there exists $g\in \ker\{\pi_1(M)\to \Z^m\}$ such that $\v(g)-\id$ is invertible over $\K$,
\en
 then $\norm_{\tau(M,\v)}$ is a seminorm on $\hom(\R^m,\R)$ and for any $\phi:\R^m\to \R$ we
 have
\[ ||\phi\circ \psi||_T \geq ||\phi||_{\tau(M,\v)}. \]
\end{theorem}

We point out that if $g\in \ker\{\pi_1(M)\to \Z^m\}$, then
$\v(g)-\id$ is defined over $\K$ since $\v$ is $\psi$--compatible.
We refer to $\norm_{\tau(M,\v)}$ as the \emph{higher--order
Alexander norm}.

In the case that $\ktbig$ equals $\Q\tbigpm$, the usual commutative Laurent polynomial ring, we
recover McMullen's Alexander norm $\norm_A$ (cf. \cite{Mc02}).
 The general commutative case is the
main result in \cite{FK05b}. The proof we give here is different in its nature from the proofs
in \cite{Mc02} and  \cite{FK05b}.

\begin{proof}
 In the case that
$m=1$ it is clear that $\norm_{\tau(M,\v)}$ is a seminorm. The fact that it gives a lower bound
on the Thurston norm was shown in \cite{Co04,Ha05,Tu02b,Fr05}. We therefore assume now that
$m>1$.

We first  show that $||\phi\circ \psi||_T \geq ||\phi||_{\tau(M,\v)}$ for any $\phi:\R^m\to \R$.
 Since both sides are $\N$--linear and continuous we only have to show that
$||\phi\circ \psi ||_T \geq ||\phi||_{\tau(M,\v)}$ for all epimorphisms $\phi:\Z^m\to \Z$. So
let $\phi:\Z^m\to \Z$ be an epimorphism.

Pick $\mu\in \Z^m$ with $\phi(\mu)=1$ as in the definition of $\deg_{\phi}(\tau(M,\v))$. We can
then again form the rings $\K[\ker(\phi)]\spm$ and $\K(\ker(\phi))(s)$. First note that by
Theorem \ref{thm:kstdet2}
\[ ||\phi||_{\tau(M,\v)}=\deg_{\phi}(\tau(M,\v))  \]
since $\phi$ is surjective. The representation
\[ \pi_1(M)\to
\gl(\ktbig,d)\to \gl(\k(\ker(\phi))\spm,d)\] is $\phi$--compatible
since $\pi_1(M)\to \gl(\ktbig,d)$ is $\psi$--compatible.
 It now follows from \cite[Theorem~1.2]{Fr05} that $||\phi\circ \psi||_T \geq
\deg(\tau(M,\K(\ker(\phi))(s)))=\deg_{\phi}(\tau(M,\v))$ (cf. also \cite{Tu02b}).

In the remainder of the proof we will show that if $m>1$ then the
Reidemeister torsion $\tau(M,\v)\in K_1(\ktbigfield)/\pm
\v(\pi_1(M))$ can be represented by a matrix defined over $\ktbig$.
It then follows from Theorem \ref{thm:kstdet} that
$\norm_{\tau(M,\v)}$ is a seminorm.

First consider the case that $\v$ is a commutative representation.  Let $\F$ be the commutative
subfield $\F$ in the definition of a commutative representation. Denote by $\F\tbigpm$ the
ordinary Laurent polynomial ring.
 Then we have $\psi$--compatible representations
$\pi_1(M)\to \gl(\F\tbigpm,d)\hookrightarrow \gl(\K\tbigpm,d)$. By \cite[Proposition~3.6]{Tu01}
we have
\[ \tau(M, \F(t_1,\dots,t_m))=\tau(M, \ktbigfield) \in K_1(\ktbigfield)/\pm \v(\pi_1(M)).\]
Since $m>1$ it follows from \cite[Theorem~4.7]{Tu01} combined with \cite[Lemmas~6.2 and
6.5]{FK05b} that $\det(\tau(M, \F(t_1,\dots,t_m)))\in \F(t_1,\dots,t_m)$ equals the twisted
multivariable Alexander polynomial, in particular it is defined over $\F[t_1^{\pm
1},\dots,t_m^{\pm 1}]$. This concludes the proof in the commutative case.

It therefore remains to consider the case that there exists $g\in \ker\{G\to \Z^m\}$ such that
$\v(g)-\id$ is invertible.
 We first consider the case that $M$ is a closed
3--manifold. Let $h=g$. Now pick a Heegard decomposition  $M=G_0\cup H_0$. We can add a handle
to ${G_0}$ in $M\sm {G_0}$ so that the core represents $g$. Adding further handles in $M\sm
{G_0}$ we can assume that the complement is again a handlebody. We call the two handlebodies
$G_1$ and $H_1$.

Now we can add a handle to ${H_1}$ in $M\sm {G_1}$ so that the core represents $h$. Adding
further handles in $M\sm {H_1}$ we can assume that the complement is again a handlebody. We call
the two handlebodies $G$ and $H$. Note that $g$ is still represented by a handle of $G$. Now
give $M$ the CW structure as follows: Take one 0--cell, attach $1$--cells along a choice of
cores of $G$ such that $g$ corresponds to one $1$--cell. Attach 2--cells along cocores of $H$
such that one cocore corresponds to $h$. Finally attach one  3--cell.

Denote the number of 1--cells by $n$. Consider the chain complex of the universal cover
$\ti{M}$:
\[ 0\to C_3(\ti{M})^1\xrightarrow{\partial_3} C_2(\ti{M})^n\xrightarrow{\partial_2}
C_1(\ti{M})^n\xrightarrow{\partial_1} C_0(\ti{M})^1\to 0,\]
 where the supscript indicates the
rank over $\Z[\pi_1(M)]$. Picking appropriate lifts of the cells of $M$ to cells of $\ti{M}$ and
picking an appropriate order we get bases  for the $\Z[\pi_1(M)]$--modules $C_i(\ti{M})$, such
that if $A_i$ denotes the matrix corresponding to $\partial_i$, then $A_1$ and $A_3$ are of the
form
\[ \ba{rcl} A_3&=&(1-g,1-g_2,\dots,1-g_n)^t, \\
  A_1&=&(1-h,1-h_2,\dots,1-h_n), \ea \]
for some $g_i,h_i \in \pi_1(M), i=2,\dots,n$. By assumption $\id-\v(g)$ and $\id-\v(h)$ are
invertible over $\K$.  Denote by $B_2$ the result of deleting the first column and the first row
of $A_2$. Let $\tau:=(\id-\v(g))^{-1}\v(B_2)(\id-\v(h))^{-1}$. Note that $\tau$ is defined over
$\ktbig$. Since we assume that $\tau(M,\v)\ne 0$ it follows that $\v(B_2)$ is invertible over
$\ktbigfield$ and $\tau(M,\v)=\tau \in K_1(\ktbigfield)/\pm \v(\pi_1(M))$ (we refer to
\cite[Theorem~2.2]{Tu01}   for details). Therefore $\tau(M,\v)\in K_1(\ktbigfield)/\pm
\v(\pi_1(M))$ can be represented by a matrix defined over $\ktbig$.

In the case that $M$ is a 3--manifold with non--empty toroidal boundary we can find a (simple)
homotopy equivalence to a 2--complex $X$ with $\chi(X)=0$. We can assume that the CW--structure
has one 0--cell, $n$ 1--cells and $n-1$ 2--cells, furthermore we can assume that one of the
1--cells represents an element $h\in \ker\{\psi:G\to \Z^m\}$ such that $\id-\v(h)$ is
invertible. We get a chain complex
\[ 0\to  C_2(\ti{X})^{n-1}\xrightarrow{\partial_2}
C_1(\ti{X})^n\xrightarrow{\partial_1} C_0(\ti{X})^1\to 0.\] Picking appropriate lifts of the
cells of $X$ to cells of $\ti{X}$ we get bases
 for the $\Z[\pi_1(X)]$--modules
$C_i(\ti{X})$, such that if $A_i$ denotes the matrix corresponding to $\partial_i$, then $A_1$
is of the form
\[ \ba{rcl}
  A_1&=&(1-h,1-h_2,\dots,1-h_n), h_i \in \pi_1(M). \ea \]
Now denote by $B_2$ the result of deleting the first row of $A_2$. Then $\tau:=
\v(B_2)(\id-\v(h))^{-1}$ is again defined over $\ktbig$ and the proof continues as in the case
of a closed 3--manifold.
\end{proof}

\begin{remark}
Note that if follows from \cite{Fr05} that if $M$ is closed, or if
$M$ has toroidal boundary, then $\tau(M,\v)\ne 0$ is equivalent to
$H_1(M;\ktbigfield)=0$, or equivalently, that $H_1(M;\ktbig)$ has
rank zero over $\ktbig$.
\end{remark}

\begin{remark}
Note that the computation of $f_d\in \ktbig$ and $f_n\in \ktbig$
such that $\det(\tau(M,\v))=f_nf_d^{-1}$ is computationally
equivalent to the computation of $\deg_\phi(\tau(M,\v))$ for some
$\phi:H_1(M)\to \Z$. Put differently we get the perhaps surprising
fact that computing the higher--order Alexander norm  does not take
longer than computing  a single higher--order one--variable
Alexander polynomial.
\end{remark}

\section{Examples of $\psi$--compatible homomorphisms}

\subsection{Skew fields of group rings} \label{section:skewfield}
A group $G$ is called locally indicable if for every finitely generated subgroup $U\subset G$
there exists a non--trivial homomorphism $U\to \Z$.

\begin{theorem}\label{thm:tfa}
Let $G$ be a locally indicable and amenable group and let $R$ be a subring of $\C$. Then $R[G]$
is an Ore domain, in particular it embeds in its classical right ring of quotients $\K(G)$.
\end{theorem}

It follows from \cite{Hi40}  that $R[G]$ has no zero divisors. The theorem  now follows from
\cite{Ta57} or \cite[Corollary~6.3]{DLMSY03}.

A group $G$ is called poly--torsion--free--abelian (PTFA) if there exists a filtration
\[ 1=G_0 \subset G_1\subset \dots \subset G_{n-1}\subset G_n=G \]
such that $G_{i}/G_{i-1}$ is torsion free abelian. It is well--known that PTFA groups are
amenable and locally indicable (cf. \cite{St74}). The group rings of PTFA groups played an
important role in \cite{COT03}, \cite{Co04} and \cite{Ha05}.


\subsection{Admissible pairs and multivariable skew Laurent polynomial rings}
We slightly generalize a definition from \cite{Ha06}.

\begin{definition} Let $\pi$ be a group and let $\psi:\pi \to \Z^m$ be an epimorphism and let
$\varphi:\pi\to G$ be an epimorphism to a locally indicable and amenable group $G$ such that
there exists a map $G\to \Z^m$ (which we also denote by $\psi$) such that
 \[   \xymatrix {
 \pi \ar[dr]_{\psi} \ar[r]^{\varphi} &G \ar[d]^{\psi}\\& \Z^m }
     \]
 commutes. Following \cite[Definition~1.4]{Ha06} we call $(\varphi,\psi)$
an {\em admissible pair} for $\pi$.
\end{definition}

Clearly $G_\psi:=\ker\{G\to \Z^m\}$ is locally indicable and amenable. It follows now from
\cite[Lemma~3.5 (ii),
 p.~609]{Pa85} that $(\Z[G],\Z[G_\psi]\sm \{0\})$ satisfies the Ore property.
Now pick elements $t^\a \in G, \a \in \Z^m$ such that $\psi(t^\a)=\a$ and $t^{n\a}=(t^\a)^{n}$ for
any $\a\in \Z^m,n\in \Z$.

 Clearly  $\Z[G](\Z[G_\psi]\sm \{0\})^{-1}=\sum_{\a\in \Z^m}\K(G_\psi)t^\a$ is a
multivariable skew Laurent polynomial ring of rank $m$ over the
field $\K(G_{\psi})$ as defined in Section \ref{section:multivar}.
We denote this ring by $\K(G_\psi)\tbigpm$. Note that $\Z[\pi]\to
\Z[G]\to \K(G_\psi)\tbigpm$ is a $\psi$--compatible homomorphism and
that  $\K(G_\psi)\tbigfield$ is canonically isomorphic to $\K(G)$.

A family of examples of admissible pairs is provided by the rational derived series of a group
$\pi$ introduced by the second author (cf. \cite[Section~3]{Ha05}). Let $\pi_r^{(0)}:=\pi$ and
define inductively
\[ \pi_r^{(n)}:=\big\{ g\in \pi_r^{(n-1)} | \, g^d \in \big[\pi_r^{(n-1)},\pi_r^{(n-1)}\big] \mbox{ for some }d\in \Z \sm \{0\} \big\}.\]
Note that $\pi_r^{(n-1)}/\pi_r^{(n)}\cong
\big(\pi_r^{(n-1)}/\big[\pi_r^{(n-1)},\pi_r^{(n-1)}\big]\big)/\mbox{$\Z$--torsion}$.
 By \cite[Corollary~3.6]{Ha05} the quotients $\pi/\pi_r^{(n)}$ are PTFA groups
for any $\pi$ and any $n$. If $\psi:\pi\to \Z^m$ is an epimorphism,
then  $(\pi\to \pi/\pi_r^{(n)},\psi)$ is an admissible pair for
$\pi)$ for any $n>0$.

\subsection{Admissible pairs and seminorms} \label{section:rhon}
Let $M$ be a 3--manifold with empty or toroidal boundary.
Let  $(\v:\pi_1(M)\to G,\psi:\pi_1(M)\to \Z^m)$  be an admissible pair for $\pi_1(M)$. We denote
the induced map $\Z[\pi_1(M)]\to \K(G_\psi)\tbigfield$ by $\v$ as well.

 Let $\phi:\Z^m\to \Z$ be a non--trivial homomorphism. We
denote the induced homomorphism $G\to \Z^m\to \Z$ by $\phi$ as well. We write $G_\phi:=\ker\{G\to
\Z\}$. Pick $\mu\in G$ such that $\phi(\mu)\Z=\im(\phi)$. We  define $\Z[G_\phi]\upm$ via $uf=\mu
f\mu^{-1}u$. Note that we get an isomorphism $\K(G_\phi)(u)\cong \K(G)$. If $\tau(M,\v)\ne 0$, then
we define
\[ \ol{\delta}_{G}(\phi):=\max\{0,\deg(\tau(M,\K(G_\phi)(u)))\} \]
otherwise we write $\ol{\delta}_{G}(\phi)=-\infty$. We will adopt
the convention that $-\infty<a$ for any $a\in \Z$. By \cite{Fr05}
this agrees with the definition in \cite[Definition~1.6]{Ha06} if
$\delta_G(\phi)\ne -\infty$ and if  $\v:G\to \Z^m$ is not an
isomorphism or if $m>1$.  In the case that $\v:G\to \Z$ is an
isomorphism and $M\ne S^1\times D^2, S^1\times S^2$, this definition
differs from \cite[Definition~1.6]{Ha06} by the term $1+b_3(M)$. In
the case that $\v:\pi\to \pi/\pi_r^{(n+1)}$ then we also write
$\ol{\delta}_{n}(\phi)=\ol{\delta}_{\pi/\pi_r^{(n+1)}}(\phi)$.

\begin{theorem} \label{thm:harveynorm}
Let $M$ be a 3--manifold with empty or toroidal boundary.
Let $(\v:\pi_1(M)\to G,\psi:\pi_1(M)\to \Z^m)$ be an admissible pair for $\pi_1(M)$
such that $\tau(M,\v)\ne 0$.
 Then for any $\phi:\Z^m\to \Z$ we have
\[ ||\phi||_{\tau(M,\v)}=\ol{\delta}_G(\phi), \]
and $\phi\mapsto \max\{0,\ol{\delta}_G(\phi)\}$ defines a seminorm which is a lower bound on the
Thurston norm.
\end{theorem}

Note that this theorem implies in particular Theorem \ref{thm:mainintro}.

\begin{proof}
Let $\phi:\Z^m\to \Z$ be a non--trivial homomorphism. As in Section \ref{section:multivar} we
can form $\K(G_\phi)\spm$ and $\K(G_\psi)(\ker(\phi))\spm$. Note that these rings are
canonically isomorphic Laurent polynomial rings. If $\psi:G\to \Z^m$ is an isomorphism, then
$\v$ is commutative. Otherwise we can find a non--trivial $g\in \ker(\psi)$, so clearly
$1-\v(g)\ne 0\in \K(G)$. This shows that we can apply Theorem \ref{mainthmnorm} which then
concludes the proof.
\end{proof}

In the case that $\v:\pi\to \pi/\pi_r^{(n+1)}$ we denote the
seminorm $\phi\mapsto \max\{0,\ol{\delta}_{n}(\phi)\}$ by $\norm_n$.
Note that in the case $n=0$ this was shown by the second author
\cite[Proposition~5.12]{Ha05} to be equal to McMullen's Alexander
norm \cite{Mc02}.

\subsection{Admissible triple} \label{section:admissible}

We now slightly extend  a definition from \cite{Ha06}.

\begin{definition}
Let $\pi$ be a group and $\psi:\pi\to \Z^m$ an epimorphism. Furthermore let $\varphi_1:\pi \to
G_1$ and $\varphi_2:\pi \to G_2$ be epimorphisms to  locally indicable and amenable groups $G_1$
and $G_2$. We call $(\varphi_1,\varphi_2,\psi)$ an {\em admissible triple} for $\pi$ if there
exist epimorphisms $\Phi:G_1\to G_2$ and $\psi_2:G_2\to \Z^m$ such that $\varphi_2=\Phi\circ
\varphi_1$, and $\psi=\psi_2\circ \varphi_2$.
\end{definition}

Note that in particular $(\varphi_i,\psi), i=1,2$ are admissible pairs for $\pi$. Combining
Theorem \ref{thm:harveynorm} with \cite[Theorem~1.3]{Fr05} (cf. also \cite{Ha06}) we get the
following result.

\begin{theorem}\label{thm:triple}
Let $M$ be a 3--manifold with empty or toroidal boundary.
If $(\varphi_1,\varphi_2,\psi)$ is an admissible triple for $\pi_1(M)$  such that
$\tau(M,\v_2)\ne 0$, then we have the following inequalities of seminorms:
\[ \norm_{\tau(M,\v_2)}\leq \norm_{\tau(M,\v_1)} \leq \norm_T. \]
In particular we have
\[ \norm_0\leq \norm_1 \leq \dots\leq \norm_T.\]
\end{theorem}

Let $M$ be a 3--manifold with empty or toroidal boundary and let $\phi\in H^1(M;\Z)$. Since $\delta_n(\phi)\in \N$ for all $n$ it follows immediately from Theorem \ref{thm:triple} that there 
exists $N\in \N$ such that $\delta_n(\phi)=\delta_N(\phi)$ for all $n\geq N$. But we can in fact prove a slightly stronger statement, namely that there
exists such an $N$ independent of the choice of $\phi \in H^1(M;\Z)$. 

\begin{proposition} \label{prop:eventconstant}
Let $M$ be a 3--manifold with empty or toroidal boundary.
There
exists  $N\in \N$ such that $\delta_n(\phi)=\delta_N(\phi)$ for all $n\geq N$ and all $\phi \in H^1(M;\R)$.
\end{proposition}

\begin{proof}
Write $\pi=\pi_1(M), \pi_n=\pi/\pi_r^{(n+1)}$ and $m=b_1(M)$. Let $\psi:\pi\to \Z^m$ be an epimorphism. Write $(\pi_n)_{\psi}=\ker\{\psi:\pi_n\to \Z^m\}$.
Now pick elements $t^\a \in \pi_n, \a \in \Z^m$ such that $\psi(t^\a)=\a$ and $t^{k\a}=(t^\a)^{k}$ for
any $\a\in \Z^m,k\in \Z$.
Consider the map $\Z[\pi]\to
\Z[\pi_n]\to \K((\pi_n)_\psi)\tbigfield$. We write
$\tau_n=\tau(M,\K((\pi_n)_\psi)\tbigfield)$. We can find $f_n,g_n\in \K((\pi_n)_\psi)\in \tbigpm$ such that
$\tau_n=f_ng_n^{-1}$. 

Given a seminorm $s$ on $H^1(N;\R)$ whose normball is a (possibly non--compact) polygon we can study its dual polytope
$d(s)$. Note that given $f=\sum_{\a\in \Z^m}a_\a t^\a \in \K((\pi_n)_\psi)\in \tbigpm$ the dual polytope $d(\norm_f)$ equals the Newton polygon $N(f)$
which is the convex hull of $\{ \a | a_\a \ne 0\}$.  Clearly $d(\norm_f)$ has only integral vertices.

By the definition of $\delta_n=\norm_{\tau_n}=\norm_{f_gg_n^{-1}}$ it follows that
\[ d(\delta_n)+d(g_n)=d(\tau_n)+d(g_n)=d(f_n)\]
where $``+''$ denotes the Minkowski sum of convex sets. It is easy to see that this implies that $d(\delta_n)$ has only integral vertices.

Theorem \ref{thm:triple} implies that there is a sequence of inclusions
\[ d(\delta_0) \subset d(\delta_1)\subset \dots \subset d(\norm_T).\]
Since $d(\norm_T)$ is compact and since $d(\delta_n)$ has integral vertices for all $n$ it follows immediately that
there
exists  $N\in \N$ such that $d(\delta_n)=d(\delta_N)$ for all $n\geq N$.  This completes the proof of the proposition.
\end{proof}

\section{Examples} \label{sec:example}

Before we discuss the Thurston norm of a family of links we first need to introduce some
notation for knots. Let $K$ be a knot. We denote the knot complement by $X(K)$. Let
 $\phi:H_1(X(K))\to \Z$ be an isomorphism. We write
$\ol{\delta}_n(K):=\ol{\delta}_n(\phi)$. This agree with the
original definition of Cochran \cite{Co04} for $n>0$ and if
$\Delta_K(t)=1$, and it is one less than Cochran's definition
otherwise.

\begin{figure}[h]
\begin{center}
\begin{picture}(250,130)
\put(0,-10){\includegraphics[scale=0.7]{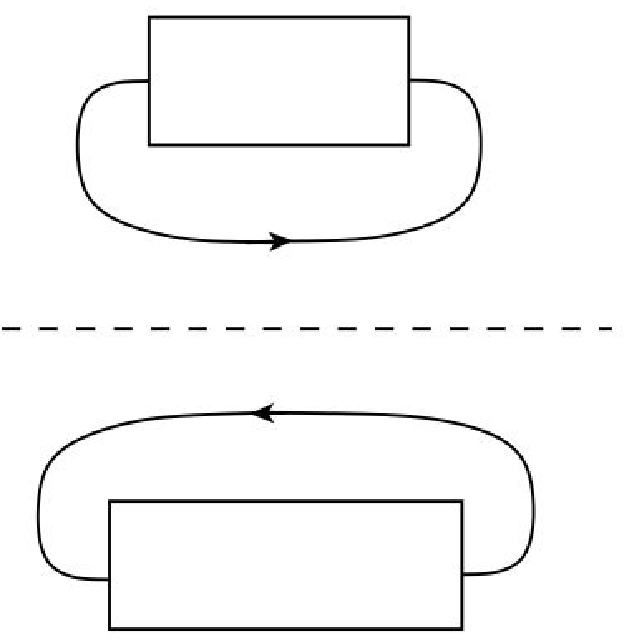}}
\put(55,12){$L$}
\put(55,110){$K$}
\put(150,0){\includegraphics[scale=0.7]{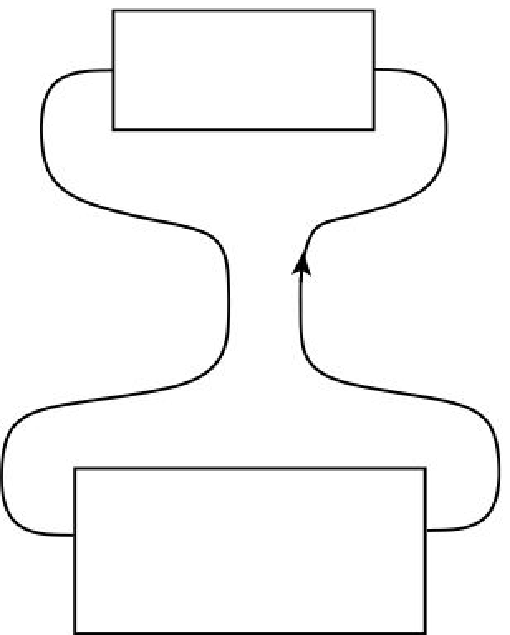}}
\put(195,12){$L$}
\put(195,110){$K$}
\end{picture}
 \caption{The link $L\#_i K$.} \label{hopflink}
\end{center}
\end{figure}

In the following let $L=L_1\cup\dots\cup L_m$ be any ordered
oriented $m$--component link. Let $i\in \{1,\dots,m\}$. Let $K$ be
an oriented  knot with $\Delta_K(t)\ne 1$ which is separated from
$L$ by a sphere $S$. We pick a path from a point on $K$ to a point
on $L_i$ and denote by $L\#_i K$ the link given by performing the
connected sum of $L_i$ with $K$ (cf. Figure \ref{hopflink}). Note
that this connected sum is well--defined, i.e. independent of the
choice of the path. We will study the Thurston norm of $L\#_i K$.

Now assume that $L$ is a non--split link with at least two
components and such that $\norm_0=\norm_T$. Many examples of such
links are known (cf. \cite{Mc02}). For the link $L\#_i K$ denote its
meridians by $\mu_i,i=1,\dots,m$. Let $\psi:H_1(X(L\#_i K))\to \Z^m$
be the isomorphism given by $\psi(\mu_i)=e_i$, where $e_i$ is the
$i$--th vector of the standard basis of $\Z^m$.

We write $\pi:=\pi_1(X(L\#_i K))$. For all $\a\in \Z^m$ we pick
$t^\a \in \pi/\pi_r^{(n+1)}$ with $\psi(t^\a)=\a$ and such that
$t^{l\a}=(t^\a)^l$ for all $\a\in \Z^m$ and $l\in \Z$. Furthermore
write $t_i:=t^{e_i}$.

\begin{proposition}\label{propalexlink}
Consider the natural map $$\v:\pi\to
\K(\pi/\pi_r^{(n+1)})=\k(\pi_\psi/\pi_r^{(n+1)})(t_1,\dots,t_m).$$
where $\pi$ is as defined above.  There exists an element $f(t_i)\in \K(\pi_{\psi}/\pi_r^{(n+1)})[t_i^{\pm
1}]\subset \k(\pi_\psi/\pi_r^{(n+1)})[t_1^{\pm 1},\dots,t_m^{\pm
1}]$ such that $\deg(f(t_i))=\ol{\delta}_n(K)+1$, and there exists a
$d=d(t_1,\dots,t_m)\in \ktbigfield$ with $\norm_d=\norm_{0}$, such
that \be \label{equ:comptau} \tau(X(L\#_i
K),\v)=d(t_1,\dots,t_m) f(t_i) \in
K_1(\k(\pi_\psi/\pi_r^{(n+1)})(t_1,\dots,t_m))/\pm \v(\pi).\ee
Furthermore, if $\ol{\delta}_n(K)=2\gen(K)-1$, then
\[ \norm_{\tau(X(L\#_i K),\v)}=\norm_T.\]
\end{proposition}


\begin{proof}

Let $S$  be the embedded sphere in $S^3$ coming from the definition
of the connected sum operation (cf. Figure \ref{hopflink}). Let $D$
be the annulus $S\cap X(L\#_i K)$ and we denote by $P$ the closure
of the component of $X(L\#_i K)\sm D$ corresponding to $K$. We
denote the closure of the other component by $P'$ (see Figure
\ref{linkcut} below). Note that $P$ is homeomorphic to $X(K)$ and
$P'$ is homeomorphic to $X(L)$.
\begin{figure}[h]
\begin{center}
\begin{picture}(100,120)
\put(0,0){\includegraphics[scale=0.7]{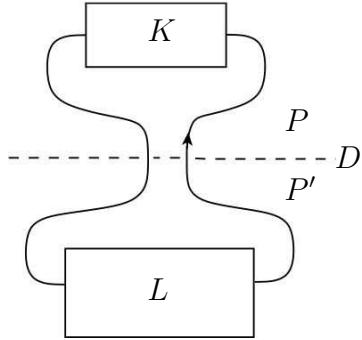}}
\put(55,14){$L$}
\put(55,112){$K$}
\put(107,52){$P^{\prime}$}
\put(107,78){$P$}
\put(126,64){$D$}
\end{picture}
 \caption{The  link complement of $L\#_i K$
cut along the annulus $D$.} \label{linkcut}
\end{center}
\end{figure}
Denote the induced maps to
$\mathcal(K):=\K(\pi_\psi/\pi_r^{(n+1)})(t_1,\dots,t_m)$ by $\v$ as
well. We get an exact sequence
\[
0 \to  C_*^\v(D;\mathcal(K)) \to C_*^\v(P;\mathcal(K))\oplus
C_*^\v(P';\mathcal(K)) \to C_*^\v(X(L\#_i K);\mathcal(K)) \to 0 \]
of chain complexes. It follows from \cite[Theorem~3.4]{Tu01} that
\be \label{equ:tau} \tau(P,\v)\tau(P',\v)= \tau(D,\v)\tau(X(L_i\#
K),\v)\in \big(K_1(\mathcal(K))/\pm \v(\pi)\big) \,\cup \{0\}.\ee
First note that $D$ is homotopy equivalent to a circle  and that
$\im\{\psi:\pi_1(D)\to \Z^m\}=\Z e_i$. It is now easy to see that
$\tau(D,\v)=(1-at_i)^{-1}$ for some $a\in
\k(\pi_\psi/\pi_r^{(n+1)})\sm \{0\}$.

Next note that $\im\{\psi:\pi_1(P)\to \Z^m\}=\Z e_i$. In particular $\tau(P,\v)$ is defined over
the one--variable Laurent polynomial ring $\K(\pi_\psi/\pi_r^{(n+1)})[t_i^{\pm 1}]$ which is a
PID. Recall that we can therefore assume that its Dieudonn\'e determinant $f(t_i)$ lies in
$\K(\pi_\psi/\pi_r^{(n+1)})[t_i^{\pm 1}]$ as well.
\begin{claim}
\[ \deg(\tau(P,\v:\pi_1(P)\to \K(\pi_\psi/\pi_r^{(n+1)})(t_i))=\ol{\delta}_n(K).\]
\end{claim}
First recall that there exists a homeomorphism $P \cong X(K)$. We
also have an inclusion $X(L\#_i K)\to X(L_i\# K)$. Combining with
the degree one map $X(L_i\# K)\to X(K)$ we get a factorization of an
automorphism of $\pi_1(X(K))$ as follows:
\[ \pi_1(X(K))\cong \pi_1(P)\to \pi_1(X(L\#_i K))\to \pi_1(X(L_i\# K))\to \pi_1(X(K)).\]
Since the rational derived series is functorial (cf. \cite{Ha05}) we in fact get that
\[ \ba{rcl} \pi_1(X(K))/\pi_1(X(K))_r^{(n+1)}&\cong &\pi_1(P)/\pi_1(P)_r^{(n+1)}\\
&\to &\pi_1(X(L_i\# K))/\pi_1(X(L_i\# K))_r^{(n+1)}\\&\to&
\pi_1(X(K))/\pi_1(X(K))_r^{(n+1)}\ea
\] is an isomorphism.
 In particular \[
\pi_1(X(K))/\pi_1(X(K))_r^{(n+1)} \to \pi_1(X(L\#_i
K))/\pi_1(X(L\#_i K))_r^{(n+1)}\] is injective, and the induced map
on Ore localizations is injective as well. Finally note that
$\ker\{\pi_1(X(K))\to \pi_1(P)\xrightarrow{\psi} \Z^m\}=\ker(\phi)$
where $\phi:\pi_1(X(K))\to \Z$ is the abelianization map. It now
follows that
\[ \ba{rcl}
\ol{\delta}_n(K)
&=& \deg(\tau(X(K),\pi_1(X(K))\to \K(\pi_1(X(K))_\phi/\pi_1(X(K))_r^{(n+1)})(t_i))\\
&=&
\deg(\tau(X(K),\pi_1(X(K))\to \K(\pi_\psi/\pi_r^{(n+1)})(t_i))\\
&=&\deg(\tau(P,\pi_1(P)\to \K(\pi_\psi/\pi_r^{(n+1)})(t_i)).
\ea \]
Note that the second equality follows from the functoriality of torsion (cf.
\cite[Proposition~3.6]{Tu01}) and the fact that going to a supfield does not change the degree
of a rational function. This concludes the proof of the claim.

\begin{claim}
We have the following equality of norms on $H^1(X(L);\Z)$:
\[ \norm_{\tau(P',\v)}=\norm_{T}.\]
\end{claim}

First recall that $P'$ is homeomorphic to $X(L)$.  The claim now follows immediately from Theorem
\ref{thm:triple} applied to $\v$ and to the abelianization map of $\pi_1(P')$, and from the
assumption that $\norm_{0}=\norm_{T}$ on $H^1(X(L);\Z)$.

Putting these computations together and using Equation
(\ref{equ:tau}) we now get a proof of Equation (\ref{equ:comptau}).

Now assume that  $\ol{\delta}_n(K)=2\gen(K)-1$. Let $S_i$ be a
Seifert surface of $K$ with minimal genus. Let $\phi:\Z^m\to \Z$ be
an epimorphism and let $l=\phi(\mu_i)\in \Z$. We first view $\phi$
as an element in $\hom(H_1(X(L);\Z)$. A standard argument shows that
$\phi$ is dual to a (possibly disconnected) surface $S$ which
intersects the tubular neighborhood of $L_i$ in exactly $l$ disjoint
curves. Then the connected sum $S'$ of $S$ with $l$ copies of $S_i$
gives a surface in $X(L\#_i K)$ which is dual to $\phi$ viewed as an
element in $\hom(H_1(X(L\#_i K);\Z)$. A standard argument shows that
$S'$ is Thurston norm minimizing (cf. e.g. \cite[p.~18]{Lic97}).

Clearly $\chi(S')=\chi(S)+l(\chi(S_i)-1)$. A straightforward
argument shows that furthermore
$\chi_-(S')=\chi_-(S)+l(\chi_-(S_i)+1)$ since $L$ is not a split
link and since $K$ is non--trivial.

We now compute
\[ \ba{rcl} ||\phi||_T&=&\chi_-(S')\\
&=&\chi_-(S)-n(\chi(S_i)-1)\\
&=&||\phi||_{T}+2l\gen(K)\\
&=&||\phi||_{d}+2(\ol{\delta}_n(K)+1)\\
&=&||\phi||_d+2\deg(f(t_i))\\
&=&||\phi||_{\tau(X(L\#_i K),\v)}.\ea\] By the $\R$--linearity and the continuity of the norms
it follows that
\[ ||\phi||_{\tau(X(L\#_i K),\v)}=||\phi||_T\]
for all $\phi:\Z^m\to \R$.
\end{proof}

Denote by $\lozenge(n,m)$ the convex polytope given by the vertices $(\pm \frac{1}{n},0)$ and
$(0,\pm \frac{1}{m})$. Let $(n_i)_{i\in \N}$ and $(m_i)_{i\in \N}$ be never decreasing sequences
of odd positive  numbers which are eventually constant, i.e. there exists an $N$ such that
$n_i=n_N$ for all $i\geq N$ and $m_i=m_N$ for all $i\geq N$. According to \cite{Co04} we can
find knots $K_1$ and $K_2$ such that $\ol{\delta}_i(K_1)=n_i$ for any $i$,
$\ol{\delta}_N(K_1)=2\, \gen(K_1)-1$ and $\ol{\delta}_i(K_2)=m_i$ for any $i$ and
$\ol{\delta}_N(K_2)=2\, \gen(K_2)-1$.

\begin{figure}[h]
\begin{picture}(100,140)
\put(0,0){\includegraphics[scale=0.7]{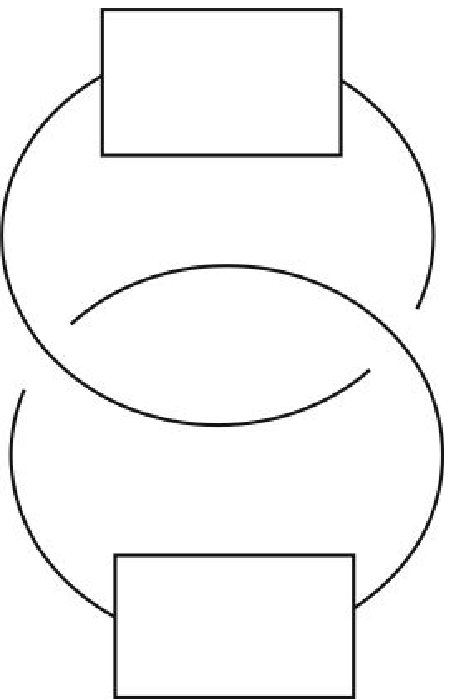}}
\put(40,120){$K_{1}$}
\put(40,12){$K_{2}$}
\end{picture}
 \caption{$H(K_{1},K_{2})$ is obtained by tying $K_{1}$ and $K_{2}$ into the Hopf link} \label{fig:hopf}
\end{figure}

Let $H(K_1,K_2)$ be the link formed by adding the two knots $K_1$
and $K_2$ from above to the Hopf link (cf. Figure \ref{fig:hopf}).
Recall that the Thurston norm ball of the Hopf link is given by
$\lozenge(1,1)$.
 Let $\pi:=\pi_1(X(L))$. It
follows immediately from applying Proposition \ref{propalexlink}
twice that the norm ball of $\norm_i$ equals $\lozenge(n_i+1,m_i+1)$
and that $\norm_N=\norm_T$. The following result is now an immediate
consequence of Proposition \ref{propalexlink}.

\begin{corollary}
We have the following sequence of inequalities of seminorms
\[ \norm_A =\norm_0\leq \norm_1 \leq \norm_2 \leq \dots \leq \norm_N =\norm_T.\]
\end{corollary}

In \cite{Ha05} the second author  gave examples of 3--manifolds $M$
such that
\[ \norm_A =\norm_0\leq \norm_1 \leq \norm_2 \leq \dots \]
but in that case it was not known whether the sequence of norms
$\norm_i$ eventually agrees with $\norm_T$.


It is an interesting question to determine which 3--manifolds satisfy $\norm_T=\norm_n$ for large enough $n$. We conclude this paper with the following
conjecture.

\begin{conjecture}
If $\pi_1(M)_{r}^{(\omega)}\equiv \bigcap_{n\in \N}\pi_1(M)_{r}^{(n)}=\{1\}$, then there exists $n\in \N$ such that $\norm_T=\norm_n$.
\end{conjecture}

\end{document}